\documentclass{article}
\usepackage{amssymb}
\usepackage{bm}
\usepackage[latin1]{inputenc} 
\usepackage{graphicx}
\usepackage[latin1]{inputenc}
\usepackage{dsfont}
\usepackage{amsmath}
\usepackage{amsthm}
\usepackage{comment}
\usepackage{natbib}
\newtheorem{theorem}{Theorem}
\newtheorem{lemma}{Lemma}
\newtheorem{corollary}{Corollary}

\newcommand{\wabs}[1]{\left|#1\right|}

\newcommand{\wrone}{\mathds R}

\newcommand{\wqab}[2]{\mathrm{Q}\!\left(#1,#2\right)}

\newcommand{\wfc}[2]{{#1}\!\left(#2\right)}

\newcommand{\wi}[1]{\wrm{i}}
\newcommand{\winter}[2]{{\mathrm{I}_{\wvec{#1},\wvec{#2}}}}
\newcommand{\winterf}[3]{\mathrm{I}_{\wvec{#1},\wvec{#2}}\!\left[{\wvec{#3}}\right]}

\newcommand{\wlr}[1]{\left( #1 \right)}
\newcommand{\wlrc}[1]{\left\{ #1 \right\}}

\newcommand{\wnorm}[1]{\left\| {#1} \right\|}

\newcommand{\wref}[1]{(\ref{#1})}
\newcommand{\wrn}[1]{{\mathds R}^{#1}}

\newcommand{\wrm}[1]{\mathrm{#1}}

\newcommand{\wrounde}[1]{\wfc{\wrm{fl}}{#1}}

\newcommand{\wset}[1]{{\left\{ #1 \right\}}}

\newcommand{\wst}[1]{\langle{#1}\rangle}

\newcommand{\wvec}[1]{\mathbf{#1}}

\begin{document}
\author{%
{\sc
Walter F. Mascarenhas\thanks{Corresponding author. Email: walter.mascarenhas@gmail.com,
supported by grant 2013/10916-2 from Fundação de Amparo à Pesquisa do Estado de São Paulo (FAPESP)}
and
André Pierro de Camargo\thanks{Email: andreuler@yahoo.com.br, supported by grant 14225012012-0 from CNPq}
}
 \\[2pt]
Institute of Mathematics and Statistics,\\ University of S\~{a}o Paulo, BRA
}
\title{On the backward stability of the second barycentric formula for interpolation}

\maketitle
\begin{abstract}
We present a new stability analysis for the second barycentric formula,
showing that this formula is backward stable when the relevant Lebesgue constant is small.
\end{abstract}

\section{Introduction}

We discuss the numerical stability of the second barycentric formula for interpolation at
nodes $x_0 < x_1 < \ldots < x_n$. This formula is given by
\begin{equation}
\label{bar2}
\wfc{q}{x;\wvec{x}, \wvec{y}, \wvec{w}} \ \ := \ \ \left.
\sum\limits_{k = 0}^{n} \frac{w_k y_k}{x - x_k} \right/ \sum\limits_{k = 0}^{n} \frac{w_k}{x - x_k}
\end{equation}
(note that we write the vector $\wvec{x} \in \wrn{n+1}$ in bold face, and $x_k$ is its $k$th coordinate.)
The function $q$ may be a polynomial in $x$ for particular choices
of the weights $w_k$, but we work in the more general context of
rational barycentric interpolation discussed in
\cite{Berrut},\cite{Bos}, \cite{Bos2}, \cite{Floater} and \cite{Hormann}.

In practice, we approximate a function $f: [x^-,x^+] \rightarrow \wrone{}$ using
the formula $q$ in \wref{bar2} in three steps:
\begin{itemize}
\item[Step I:] Abstract approximation theory provides convenient nodes $x_k$ and weights $w_k$, so that,
in exact arithmetic, the error $\wfc{f}{x} - \wfc{q}{x;\wvec{x}, \wfc{f}{\wvec{x}}, \wvec{w}}$ is small for
$x \in [x^-,x^+]$.
\item[Step II:] We then obtain floating point approximations $\hat{x}_k$, $y_k$ and $\hat{w}_k$ for $x_k$, $\wfc{f}{x_k}$
and $w_k$.
\item[Step III:] Finally, we approximate $\wfc{f}{x}$ evaluating $\wfc{q}{x;\hat{\wvec{x}}, \wvec{y}, \hat{\wvec{w}}}$
numerically.
\end{itemize}

This article presents upper and lower bounds on the backward errors in Steps II and III.
We also emphasize the importance of considering the errors in these two steps.
We focus on the effects of the errors in the nodes $\wvec{x}$ and weights $\wvec{w}$,
and assume that the function values $y_k$ are exact,
because perturbations in the function values can be easily
handled using Lebesgue constants or by assuming that the perturbed
function values are of the form $y_k \wlr{1 + \beta_k}$ with $\beta_k$ small.
We discuss both the case in which the end points of the interval $[x^-,x^+]$ are nodes
and the case in which $x^-$ and $x^+$ are not nodes, and allow for the possibility
that some nodes lie outside of the interval $[x^-,x^+]$.

The overall conclusion is that formula $q$ in \wref{bar2} is backward stable when
the relevant Lebesgue constant is small, in the sense that the values
$\wfc{q}{x;\hat{\wvec{x}}, \wvec{y}, \hat{\wvec{w}}}$
obtained numerically in Steps II and III are equal to
the exact value $\wfc{q}{\tilde{x};\wvec{x}, \tilde{\wvec{y}}, \wvec{w}}$, with
$\tilde{x}$ near $x$ and $\tilde{y}_k = y_k \wlr{1 + \beta_k}$ for small $\beta_k$s.
This conclusion is different from the one presented in \cite{Hig1}, which states that
the second barycentric formula is not backward stable.
However, there is no contradiction between our conclusion and Higham's,
because we consider the favorable case in which the Lebesgue constant is small
and his conclusion applies to the worst possible scenarios.

This article has three more sections. Section \ref{section_theory}
presents upper bounds on the backward errors, showing that the second barycentric
formula is backward stable under reasonable assumptions.
Section \ref{section_lower_bound}
gives lower bounds on the backward errors, showing that, for Lagrange polynomials
and except for $\log n$ factors,
the bounds in Section \ref{section_theory} are sharp.
The last section contains a perturbation theory for the barycentric
formula, which covers the rational as well as the polynomial case. It also
presents a proof of the main theorem, which is stated in Section \ref{section_theory}.

\section{Upper bounds on the backward error}
\label{section_theory}
In this section we present upper bounds on the
backward errors in the evaluation of the second barycentric formula \wref{bar2},
complementing the bounds presented in \cite{Masc} and \cite{MascCam}.
We look at the second formula in \wref{bar2} as a linear transformation $\winter{x}{w}$
mapping $\wvec{y} \in \wrn{n+1}$ to the rational function defined by
\begin{equation}
\label{operator}
\wfc{\winterf{x}{w}{y}}{x} :=
\ \ \left\{
\begin{array}{l}
y_k \ \wrm{when} \ \ x = x_k \in \wset{x_0,x_1,\dots,x_n},   \\[0.2cm]
\left.
\sum\limits_{k = 0}^{n} \frac{w_k y_k}{x - x_k} \right/ \sum\limits_{k = 0}^{n} \frac{w_k}{x - x_k}
\ \wrm{for} \ x \in \wrone{}  \setminus \wset{x_0,\dots, x_n}.
\end{array}
\right.
\end{equation}
The letter $\mathrm{I} $ in $\winter{x}{w}$ stems from {\it Interpolant}, because
the function $\winterf{x}{w}{y}: \wrone{} \rightarrow \wrone{}$ defined by \wref{operator}
interpolates the $y_k$ at the $x_k$.

The linear map $\winter{x}{w}$ is an abstract way of looking at the second barycentric
formula \wref{bar2}, and a practical minded reader can think of $\winterf{x}{w}{y}$
as a synonym for the function $q$ in \wref{bar2}. However, by considering the linear map
$\winter{x}{w}$ we can think at a deeper level.
When the nodes and weights are such that
\begin{equation}
\label{hypo_nz_den}
\sum\limits_{k = 0}^{n} \frac{w_k}{x - x_k} \neq 0 \hspace{0.3cm}  \wrm{for} \ x \in [x^-,x^+]
\setminus \wset{x_0,\dots,x_n},
\end{equation}
the function $q$ in \wref{bar2} does not have poles in the
interval $[x^-,x^+]$, and $\winterf{x}{w}{y}$  is an element of the vector space
$\wqab{x^-}{x^+}$ of
continuous rational functions from $[x^-,x^+]$ to $\wrone{}$, in which we can define the sup norm.
It is then natural to study the norm of
$\winter{x}{w}$ with respect to the sup norm in $\wrn{n+1}$ and $\wqab{x^-}{x^+}$.
This norm is the Lebesgue constant
mentioned in the abstract, and we denote it by $\Lambda_{x^-,x^+,\wvec{x},\wvec{w}}$. Formally,
we define
\begin{equation}
\label{def_lebesgue}
\Lambda_{x^-,x^+,\wvec{x},\wvec{w}} := \wnorm{\winter{x}{w}}_\infty :=
\sup_{x \in [x^-,x^+] \ \wrm{and} \ \wvec{y} \neq 0}
\frac{\wabs{ \wfc{q}{x;\wvec{x}, \wvec{y}, \wvec{w}} }}{\wnorm{\wvec{y}}_\infty} .
\end{equation}
The articles \cite{Bos}, \cite{Bos2} and \cite{Floater} present bounds on
these Lebesgue constants, and their bounds allow us to apply the theory developed in Section \ref{section_proofs}
to the Floater-Hormann interpolants.

Throughout the article we consider a reference interval $[x^-,x^+]$, nodes
$\wvec{x}$ and weights $\wvec{w}$, and perturbed (or rounded) nodes
$\hat{\wvec{x}}$, with a corresponding interval $[\hat{x}^-,\hat{x}^+]$ and  weights $\hat{\wvec{w}}$.
Besides the Lebesgue constant, our analysis of the backward stability of the second
barycentric formula is based on the relative errors in the length of the intervals
$[x_j,x_k]$, which are measured by
\begin{equation}
\label{def_delta_jk}
\delta_{kk} := \wfc{\delta_{kk}}{\wvec{x},\hat{\wvec{x}}} := 0
\hspace{1cm} \wrm{and} \hspace{1cm}
\delta_{jk} := \wfc{\delta_{jk}}{\wvec{x},\hat{\wvec{x}}} := \frac{x_j - x_k}{\hat{x}_j - \hat{x}_k} - 1.
\end{equation}
In order to handle rounding errors in the endpoints $x^-$ and $x^+$, and
errors in nodes close to them, we also consider
\begin{equation}
\label{def_udelta}
\delta^-_j := \wfc{\delta^-_j}{x^-, \wvec{x}, \hat{x}^-, \hat{\wvec{x}}}
:= \frac{x^- - x_j}{\hat{x}^- - \hat{x}_j} - 1,
\end{equation}
\begin{equation}
\label{def_odelta}
\delta^+_j := \wfc{\delta^+_j}{\wvec{x},x^+,\hat{\wvec{x}},\hat{x}^+}
:= \frac{x^+ - x_j}{\hat{x}^+ - \hat{x}_j} - 1,
\end{equation}
with $\delta^-_j = 0$ in the particular case in which $\hat{x}^- = \hat{x}_j$,
and $\delta^+_j = 0$ when $\hat{x}^+ = \hat{x}_j$.
We combine the $\delta_j^-$, $\delta_{jk}$ and $\delta_j^+$ in the $\delta$ given by
\begin{equation}
\label{def_delta}
\delta := \max_{0 \leq j,k \leq n} \wset{\wabs{\delta^-_j}, \ \wabs{\delta_{jk}}, \  \wabs{\delta^+_j}}.
\end{equation}
Another important measure of the perturbations are the relative differences $\zeta_k$ between the
reference weights $\wvec{w}$ and the weights $\hat{\wvec{w}}$ used in computation:
\begin{equation}
\label{def_zk}
\zeta_k := \wfc{\zeta_k}{\wvec{w}, \hat{\wvec{w}}} := \frac{w_k - \hat{w}_k}{\hat{w}_k},
\end{equation}
and to avoid pathological cases we assume that
\begin{equation}
\label{first_cond}
w_k \neq 0 \hspace{1cm} \wrm{and} \hspace{1cm} \hat{w}_k \neq 0.
\end{equation}
We also make the following definitions and assumptions regarding the nodes and endpoints:
\begin{eqnarray}
x_k < x_{k+1}, \hspace{0.5cm} \\
x_k \in \wlr{x^-,x^+} \  \wrm{if \ and \ only \ if } \ \hat{x}_k \in \wlr{\hat{x}^-,\hat{x}^+}, \hspace{0.5cm} \\
x_k = x^- \ \  \wrm{if \ and \ only \ if } \ \hat{x}_k = \hat{x}^-,  \hspace{0.5cm} \\
x_k = x^+ \ \ \wrm{if \ and \ only \ if } \ \hat{x}_k = \hat{x}^+, \hspace{0.5cm} \\
k^- \ \wrm{is \ the \ smallest \ } k \ \wrm{such \ that } \ x_k > x^-, \hspace{0.5cm} \\
\label{last_cond}
k^+ \ \wrm{is \ the \ largest \ } k \ \wrm{such \ that } \ x_k < x^+, \ \ \wrm{and}\ \ k^+ \geq k^-.\hspace{0.5cm}
\end{eqnarray}

We can now state our main theorem, which provides
an upper bound on the backward errors in steps II and III.

\begin{theorem}
\label{thm_main}
Under the conditions \wref{hypo_nz_den} and \wref{first_cond}--\wref{last_cond},
let $\epsilon$ be the machine precision,
 assume that  $\wlr{2n + 5} \epsilon < 1$ and define
\begin{equation}
\label{def_main_z}
Z := \frac{\wnorm{\wfc{\bm{\zeta}}{\wvec{w},\hat{\wvec{w}}}}_\infty + \wlr{n + 2}\epsilon}{1 - \wlr{n+2}\epsilon}.
\end{equation}
If, for $\delta$ in \wref{def_delta},
\begin{equation}
\label{hypo_main_a}
\wlr{\delta + Z} \Lambda_{x^-,x^+,\wvec{x},\wvec{w}} + Z < 1,
\end{equation}
and $\hat{x} \in [\hat{x}^-,\hat{x}^+]$ is a floating point number then
the computed value $\wrounde{\wfc{q}{\hat{x};\hat{\wvec{x}}, \wvec{y},\hat{\wvec{w}}}}$ is
equal to $\wfc{q}{x;\wvec{x},\tilde{\wvec{y}},\wvec{w}}$, for some $x \in [x^-,x^+]$ such that
\begin{equation}
\label{x_near}
\wabs{x - \hat{x}} \leq \max \wset{ \wnorm{\wvec{x} - \hat{\wvec{x}}}_\infty, \
\wabs{\hat{x}^- - x^-}, \
\wabs{\hat{x}^+ - x^+}},
\end{equation}
    \begin{equation}
\label{bound_main_a}
\tilde{y}_k = y_k \wlr{1 + \alpha_k} \wlr{1 + \nu_k}
\hspace{1cm} \wrm{with} \hspace{1cm}
\wnorm{\bm{\nu}}_\infty \leq \frac{\wlr{2 n + 5} \epsilon}{1 - \wlr{2 n + 5} \epsilon}
\end{equation}
and
\begin{equation}
\label{bound_main_ba}
\wnorm{\bm{\alpha}}_\infty \leq
    \frac{\wlr{1 + \Lambda_{x^-,x^+,\wvec{x}, \wvec{w}}}\wlr{\delta + Z}}
     {1 - Z -  \wlr{\delta + Z} \Lambda_{x^-,x^+,\wvec{x},\wvec{w}}}.\\[0.2cm]
\end{equation}
\end{theorem}

Theorem \ref{thm_main} states that the value obtained by the numerical evaluation of the
second barycentric formula using approximate nodes $\hat{\wvec{x}}$ and approximate weights
$\hat{\wvec{w}}$ is the exact value corresponding to $\tilde{\wvec{y}}$ ``near'' $\wvec{y}$ and
$x$ ``near'' $\hat{x}$, according to the measures of nearness in \wref{x_near}--\wref{bound_main_ba}.
We emphasize that this theorem takes into account the fact that
both the nodes and the weights may have errors in practice, and that by disregarding
one of these errors we may underestimate the backward error.

\begin{figure}[!h]
\includegraphics[bb= -170 50 600 610, width=6.2cm, height=4.8cm]{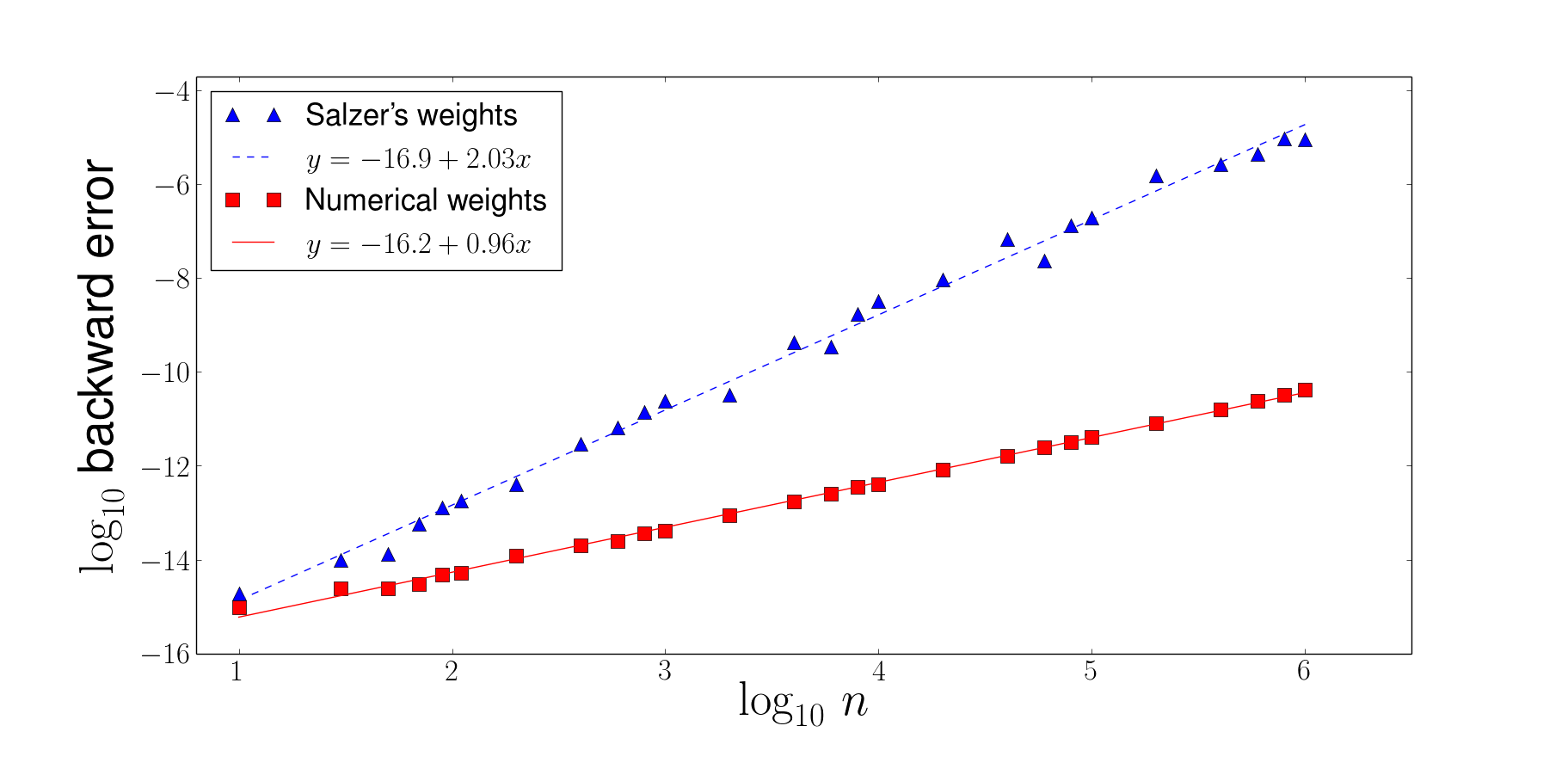}
\caption{The dependency  on the weights of the errors with rounded nodes.}
\label{figure_least_squares}
\end{figure}

Theorem \ref{thm_main} is abstract and general, and we now present examples
of its applicability in concrete situations. We state two corollaries regarding
polynomial interpolation
at the Chebyshev points of the second kind, which are defined as
\begin{equation}
\label{chebyshev_nodes}
x_k^{(c)} := - \wfc{\cos}{k \pi /n},
\end{equation}
in combination with weights obtained using the traditional formula
\begin{equation}
\label{barycentric_weights}
w_k = \wfc{\lambda_k}{\wvec{x}} := \prod_{j \neq k} \frac{1}{x_k - x_j}.
\end{equation}
We analyze two scenarios:
\begin{itemize}
\item In the first case we consider weights $\hat{\wvec{w}}$
given by the closed form expressions in \cite{SALZER}.
These weights are floating point numbers and we call them {\it Salzer's weights}.
\item In the second case we consider the weights
obtained by evaluating \wref{barycentric_weights} numerically, using the
rounded nodes $\hat{\wvec{x}}^c := \wrounde{\wvec{x}^c}$
and $\hat{\wvec{w}}  = \wrounde{\wfc{\lambda}{\hat{\wvec{x}}^c}}$,
and we call them {\it Numerical weights}.
\end{itemize}
Figure \ref{figure_least_squares} shows that these two cases are quite different
for Lagrange polynomials: although Salzer's weights contain no rounding errors,
they lead to much worse results for large $n$
(the data in this plot comes from
Tables \ref{table_salzer} and \ref{table_rounded} in Section
\ref{section_lower_bound}.) This difference is also present in Corollaries
\ref{cor_salzer_bound} and \ref{cor_numerical_bound} below, and by studying their
proof the reader will appreciate how Theorem \ref{thm_main}
can be applied in practice. Note that these corollaries provide upper bounds
on the backward error
of order $\epsilon n^2 \log n$ for Salzer's weights and $\epsilon n \log n$
for the Numerical weights,
and these numbers are in remarkable agreement with the
corresponding lines fitted by the least squares method in Figure \ref{figure_least_squares}
(recall that $\epsilon \approx 2.3 \times 10^{-16}$.)

We now present our corollaries and their proofs. In the statement of these corollaries,
$\wvec{x}^c$ are the Chebyshev nodes in \wref{chebyshev_nodes}, $\hat{\wvec{x}}^c$
are their rounded counterparts and
$\wfc{q}{x; \hat{\wvec{x}}^c,\wvec{y}, \wfc{\lambda}{\hat{\wvec{x}}^c}}$
is the $n$th degree polynomial that interpolates $\wvec{y}$ at the nodes $\hat{\wvec{x}}^c$ used in practice.

\begin{corollary}
\label{cor_salzer_bound}
If $\hat{x} \in [-1,1]$ is a floating point number, $10 \leq n \leq 2.000.000$,
$\epsilon \leq 2.3 \times 10^{-16}$,  $\wnorm{\hat{\wvec{x}}^c - \wvec{x}^c}_\infty \leq 2 \epsilon$,
and $\hat{\wvec{w}}^s$ are Salzer's weights, then there exists $x \in [-1,1]$
with $\wabs{x - \hat{x}} \leq \wnorm{\hat{\wvec{x}}^c - \wvec{x}^c}_\infty$
and $\bm{\beta} \in \wrn{n+1}$ such that
\begin{equation}
\label{bound_cor_salzer_bound_beta}
\wnorm{\bm{\beta}}_\infty \leq  3.7 \times \wlr{3 + \log n} \epsilon n^2,
\end{equation}
for which the vector $\tilde{\wvec{y}} \in \wrn{n+1}$
with entries $\tilde{y}_k = \wlr{1 + \beta_k} y_k$ satisfies
\begin{equation}
\label{bound_cor_salzer_bound}
\wrounde{\wfc{q}{\hat{x};\hat{\wvec{x}}^c, \wvec{y},\hat{\wvec{w}}^s}} =
\wfc{q}{x;\hat{\wvec{x}}^c,\tilde{\wvec{y}},\wfc{\lambda}{\hat{\wvec{x}}^c}}.\\[0.3cm]
\end{equation}
\end{corollary}

{\bf Proof of Corollary \ref{cor_salzer_bound}.}
In the context of Corollary \ref{cor_salzer_bound}, the $\delta$ in \wref{def_delta} is equal to zero,
because we consider the rounded nodes as the interpolation points from the start.
When $\wvec{w} = \wfc{\lambda}{\hat{\wvec{x}}^c}$ and $\hat{\wvec{w}} = \wfc{\lambda}{\wvec{x}^c}$,
the $\zeta_k$ in definition (17) in \cite{MascCam}
is the same as the $\zeta_k$ in definition \wref{def_zk} here
and Table 2 in that article shows that
\begin{equation}
\label{z_cor_1}
\wnorm{\bm{\zeta}}_\infty  \leq  2.4624 \wnorm{\wvec{x}^c - \hat{\wvec{x}}^c}_\infty n^2 \leq 4.9248  \epsilon n^2.
\end{equation}
Using that $10 \leq n \leq 2 \times 10^6 $, we conclude that $Z$ in  \wref{def_main_z} satisfies
\[
Z \leq \frac{4.9248 \epsilon n^2   + \wlr{n + 2} \epsilon}{1 - \wlr{2 \times 10^6 + 2} \times 2.3 \times 10^{-16}}
\leq \wlr{4.9249 n^2 + 1.0001 n + 2.0001} \epsilon
\]
\[
\leq \wlr{4.9249 + \frac{1.0001}{n} + \frac{2.0001}{n^2}}  \epsilon n^2 \leq 5.0450 \epsilon n^2
\]
\[
\leq 5.0450 \times 4 \times 10^{12} \times 2.3 \times 10^{-16}  \leq 0.0046414.
\]
Table 2 in \cite{MascCam} also shows that
\[
\Lambda_{-1,1,\hat{\wvec{x}}^c} \leq 0.67667 \log n + 1.0236
\hspace{1cm} \wrm{and} \hspace{1cm}
\Lambda_{-1,1,\hat{\wvec{x}}^c} \leq 10.841,
\]
and it follows that the $\bm{\alpha}$ in \wref{bound_main_a} satisfies
\[
\wnorm{\bm{\alpha}}_\infty \leq \frac{0.67667 \log n + 2.0236}{ 1- 0.0046414 \times 11.841} \times 5.0450  \epsilon n^2
\leq 3.6124 \wlr{3 + \log n}  \epsilon n^2.
\]
The assumptions on $\epsilon$ and $n$ lead to
\[
\wlr{2 n + 5} \epsilon \leq \wlr{4 \times 10^6 + 5} \times 2.3 \times 10^{-16} < 9.2001 \times 10^{-10},
\]
and \wref{bound_main_a} leads to $\wabs{\nu}_k \leq 1.0001 \wlr{2 n + 5} \epsilon < 9.2011 \times 10^{-10}$.
It follows that $\beta_k := \nu_k + \wlr{1 + \nu_k} \alpha_k$ satisfies
\[
\wabs{\beta_k} \leq 1.0001 \wlr{2 n + 5} \epsilon + \wlr{1 + 9.2011 \times 10^{-10}} \times  3.6124  \wlr{3 + \log n}  \epsilon n^2
\]
\[
\leq \wlr{ \frac{1.0001 \wlr{2 n + 5}}{3.6124 n^2 \wlr{3 + \log n}} + 1} \times 3.6124 \times \wlr{3 + \log n} \epsilon n^2
\leq 3.6597 \wlr{3 + \log n}  \epsilon n^2
\]
Therefore,  $\beta_k$ satisfies \wref{bound_cor_salzer_bound_beta}.
Theorem \ref{thm_main} yields $\hat{x}$ as in \wref{bound_cor_salzer_bound}
and we are done. \\[0.01cm]
\qed{}
\\[0.1cm]

\begin{corollary}
\label{cor_numerical_bound}
If $\hat{x} \in [-1,1]$ is a floating point number, $10 \leq n \leq 2.000.000$,
$\epsilon \leq 2.3 \times 10^{-16}$,  $\wnorm{\hat{\wvec{x}}^c - \wvec{x}^c}_\infty \leq 2 \epsilon$,
and $\hat{\wvec{w}}^r$ are the Numerical weights $\wrounde{\wfc{\lambda}{\hat{\wvec{x}}^c}}$,
then there exists $x \in [-1,1]$ with
$\wabs{x - \hat{x}} \leq \wnorm{\hat{\wvec{x}}^c - \wvec{x}^c}_\infty$ and
$\bm{\beta} \in \wrn{n+1}$ such that
\begin{equation}
\label{bound_cor_numerical_bound_beta}
\wnorm{\bm{\beta}}_\infty \leq
\wlr{2.2 \log n + 9.1 } \epsilon n.
\end{equation}
for which the the vector $\tilde{\wvec{y}} \in \wrn{n+1}$ with
entries $\tilde{y}_k = \wlr{1 + \beta_k} y_k$ satisfies
\begin{equation}
\label{bound_cor_numerical_bound}
\wrounde{\wfc{q}{\hat{x}; \hat{\wvec{x}}^c, \wvec{y}, \hat{\wvec{w}}^r}} =
\wfc{q}{x; \hat{\wvec{x}}^c,\tilde{\wvec{y}}, \wfc{\lambda}{\hat{\wvec{x}}^c}}.\\[0.2cm]
\end{equation}
\end{corollary}


{\bf Proof of Corollary \ref{cor_numerical_bound}.}
Lemma 3.1 in \cite{Hig1} states that
\[
w_k^{\wlr{r}} = w_k \wst{2n}_k,
\]
and using Lemma 3.1 in \cite{Hig} we conclude that the $\zeta_k$ in \wref{def_zk} satisfy
\[
\wabs{\zeta_k} = \wabs{\frac{w_k - w_k \wst{2n}_k}{w_k \wst{2n}_k}} = \frac{\wabs{1 - \wst{2n}_k}}{\wabs{\wst{2n}_k}}
= \wabs{\wst{2n}_{k'} - 1} \leq \frac{2 \epsilon n }{1 - 2 \epsilon n}
\]
\[
\leq \frac{2 \epsilon n}{1 - 2 \times 2 \times 10^6 \times 2.3 \times 10^{-16}} \leq 2.0001  \epsilon n,
\]
and the arguments after equation \wref{z_cor_1} lead to
\wref{bound_cor_numerical_bound_beta} and \wref{bound_cor_numerical_bound}.\\[0.01cm]
\qed{}

\section{Lower bounds on the backward error}
\label{section_lower_bound}

This section shows that Theorem \ref{thm_main} is sharp,
except for $\log n$ factors. These factors
are not relevant for $n$ up to one million
and we present examples in which the upper bounds provided by
Theorem \ref{thm_main} are not much larger than the maximum error
observed in practice.
We use a combination of theory and
experiments. We cannot prove that the rounding errors will be
{\bf always} above some positive number, because sometimes
the value we obtain numerically is exact. For instance,
the numerical result is exact
when we evaluate the second barycentric formula at the node $x_k$ and $y_k$ is exact,
regardless of the errors in the weights.
Therefore, we can only obtain meaningful lower bounds
under appropriate hypothesis, and experiments help us to
show that these hypothesis are fulfilled in practice.

We use Lagrange polynomials as guinea-pigs.
Since we consider reference weights $\wfc{\lambda}{\wvec{x}}$, where
$\wvec{x}$ are the nodes used in interpolation, there are no errors
in the Step I mentioned in the introduction in this case.
Moreover, we have only one $y_k$ to worry about. Formally,
Lagrange polynomials can be written in second barycentric form as
\[
\wfc{\ell_k}{x;\wvec{x}} =
\ \ \left.
 \frac{\wfc{\lambda_k}{\wvec{x}} y_k}{x - x_k} \right/ \sum\limits_{k = 0}^{n} \frac{\wfc{\lambda_k}{\wvec{x}}}{x - x_k} =
 \wfc{q}{x; \wvec{x}, \wvec{e}^k, \wfc{\lambda}{\wvec{x}}},
\]
where $y_k = 1$ and $\wvec{e}^k \in \wrn{n+1}$ is the vector with
$\wvec{e}^{(k)}_k  = 1 $  and $\wvec{e}^{(k)}_j  = 0$ for $j\neq k$.

There is a simple expression for the backward error in Step II and Step III for Lagrange
polynomials. In fact, when there is no perturbation in the nodes
and we measure the backward error in terms of the
relative perturbation in the function values,
the backward error $\beta_k$ in Steps II and III for Lagrange polynomials is such that
$\wrounde{\wfc{q}{x; \wvec{x}, \wvec{e}^k, \hat{\wvec{w}}}} = \wfc{q}{x;\wvec{x}, \wvec{e}^k \wlr{1 + \beta_k},  \wfc{\lambda}{\wvec{x}}}$,
and this condition leads to
\begin{equation}
\label{lagrange_backward}
\beta_k = \frac{\wrounde{\wfc{q}{x;\wvec{x}, \wvec{e}^k, \hat{\wvec{w}}}} -
\wfc{q}{x;\wvec{x}, \wvec{e}^k,  \wfc{\lambda}{\wvec{x}}}}{\wfc{q}{x;\wvec{x}, \wvec{e}^k,  \wfc{\lambda}{\wvec{x}}}}.
\end{equation}

This expression for $\beta_k$ allows us to prove the following theorem:
\begin{theorem}
\label{thm_lower_bound}
Assume that the $x_i$ and $\hat{w}_i$ are
floating point numbers, let $j$ and $k$ be indexes such that
$\wabs{\zeta_k} = \wnorm{\bm{\zeta}}_\infty$ and $\zeta_k \zeta_j \leq 0$ and define
\begin{equation}
\label{def_lb_s}
S = \sum_{i \neq j} \frac{\wabs{w_i}}{\wabs{x_j - x_i}}.
\end{equation}
If $2.5 (n + 3) \epsilon \leq \wnorm{\bm{\zeta}}_\infty \leq 0.001$ and $x$ is such that
\begin{equation}
\label{hypo_lb}
0 < \wabs{\frac{x - x_j}{w_j}} S \leq 0.01 \hspace{1cm} \wrm{and}\hspace{1cm}
\sup_{i \neq j} \wabs{\frac{x - x_j}{x_i - x_j}} < 0.01,
\end{equation}
then the backward error $\beta_k$ in \wref{lagrange_backward} satisfies
$\wabs{\beta_k} \geq 0.16 \, \wnorm{\bm{\zeta}}_\infty$.\\[0.5cm]
\end{theorem}

{\bf Proof of Theorem \ref{thm_lower_bound}.}
As in the proof of Theorem \ref{thm_main}, we use Stewart's relative error counter $\wst{n}$.
Equations  \wref{bar2} and  \wref{lagrange_backward} and the arguments
after equation 4.1 in \cite{Hig1} show that $\beta_k$ satisfies
\begin{equation}
\label{eq_gamma}
\ \ \frac{ \frac{ \hat{w}_k \wst{n + 3}_k}{x - x_k} }
{\sum\limits_{i = 0}^{n} \frac{\hat{w}_i \wst{n + 2}_i}{x - x_i} } \ \ = \ \
\ \ \frac{ \frac{ w_k  \wlr{1 + \beta_k} }{x - x_k} }
{\sum\limits_{i = 0}^{n} \frac{w_i}{x - x_i} }.
\end{equation}
The identities $w_i = \hat{w}_i \wlr{1 + \zeta_i}$ and $\theta_i := \wst{n+2}_i - 1$ lead to
\[
\hat{w}_i \wst{n+2}_i = w_i \frac{1 + \theta_i}{1 + \zeta_i}  = w_i \wlr{1 - \zeta_i + \psi_i},
\]
where
\[
\psi_i :=  \frac{1 + \theta_i}{1 + \zeta_i} - 1 + \zeta_i = \frac{\zeta_i^2 + \theta_i}{1 + \zeta_i}.
\]
The hypothesis $2.5 \wlr{n + 3} \epsilon \leq \wnorm{\bm{\zeta}}_\infty \leq 0.001$
and Lemma 3.1 in \cite{Hig} yield
\[
\wabs{\theta_i} \leq \frac{1}{1 - 0.001/2.5} \wlr{n+ 2} \epsilon \leq 0.401 \wnorm{\bm{\zeta}}_\infty
\]
and
\begin{equation}
\label{bound_psi}
\wabs{\psi_i} \leq \frac{0.401 + \wnorm{\bm{\zeta}}_\infty}{1 - \wnorm{\bm{\zeta}}_\infty} \wnorm{\bm{\zeta}}_\infty \leq 0.403 \wnorm{\bm{\zeta}}_\infty.
\end{equation}
Analogously, $\hat{w}_k \wst{n+3}_k = w_k \wlr{1 - \zeta_k + \phi}$
for $\phi$ such that
\begin{equation}
\label{bound_phi}
\wabs{\phi} \leq 0.403 \wnorm{\bm{\zeta}}_\infty.
\end{equation}
We can then rewrite \wref{eq_gamma} as
\[
\wlr{1 - \zeta_k + \phi} \sum\limits_{i = 0}^{n} \frac{w_i}{x - x_i}  =
\wlr{1 + \beta_k}  \sum\limits_{i = 0}^{n} \frac{w_i \wlr{1 - \zeta_i + \psi_i}}{x - x_i}
\]
and deduce that $\beta_k = N / D$ for
\begin{equation}
\label{def_xi}
\xi := \wlr{x - x_j}/w_j,
\end{equation}
\[
D := \xi \sum\limits_{i = 0}^{n} \frac{w_i \wlr{1 - \zeta_i + \psi_i}}{x - x_i}
\]
and
\[
N := \xi \wlr{1 - \zeta_k + \phi} \sum\limits_{i = 0}^{n} \frac{w_i}{x - x_i} - D.
\]
It follows that
\begin{equation}
\label{mid_d}
D = 1 - \zeta_j + \psi_j + \xi \wlr{A - B + C}
\end{equation}
for
\begin{equation}
\label{abc}
A :=  \sum_{i \neq j} \frac{w_i}{x - x_i},
\hspace{0.7cm}
B := \sum_{i \neq j} \frac{w_i \zeta_i }{x - x_i}
\hspace{0.7cm} \wrm{and} \hspace{0.7cm}
C := \sum_{i \neq j} \frac{w_i \psi_i}{x - x_i}
\end{equation}
and
\[
N = \wlr{1 - \zeta_k + \phi} \wlr{1 + \xi A} - 1 + \zeta_j - \psi_j - \xi \wlr{A - B + C}
\]
\begin{equation}
\label{nlb}
= \zeta_j - \zeta_k + \phi - \psi_j + \xi \wlr{A \phi - A \zeta_k + B - C}.
\end{equation}
The hypothesis \wref{hypo_lb} yields the bound
\[
\wabs{x - x_i} = \wabs{x_j - x_i} \wabs{1 - \frac{x_j - x}{x_j - x_i}} \geq \wabs{x_j - x_i} \wlr{1 - 0.01} \geq 0.99 \wabs{x_j - x_i},
\]
and combining this bound with \wref{hypo_lb}, \wref{def_xi} and \wref{abc} we obtain
\[
\wabs{A} \leq \frac{1}{0.99} S, \hspace{0.5cm} \wabs{B} \leq \frac{1}{0.99}\wnorm{\bm{\zeta}}_\infty S
 \hspace{0.5cm} \wrm{and} \hspace{0.5cm}
\wabs{C} \leq \frac{1}{0.99}  \wnorm{\bm{\psi}}_\infty S,
\]
for $S$ in \wref{def_lb_s}. The hypothesis \wref{hypo_lb} tells us that $\xi S \leq 0.01$ and the last equation yields
\begin{equation}
\label{bound_xia}
\wabs{\xi A} \leq \frac{0.01}{0.99}, \hspace{0.5cm} \wabs{\xi B} \leq \frac{0.01}{0.99}\wnorm{\bm{\zeta}}_\infty
 \hspace{0.5cm} \wrm{and} \hspace{0.5cm}
\wabs{\xi C} \leq \frac{0.01}{0.99}  \wnorm{\bm{\psi}}_\infty.
\end{equation}
Combining this estimate with the hypothesis that $\zeta_k$ and $\zeta_j$ have opposite signs and using
\wref{nlb} and reminding that $\wabs{\zeta_k} = \wnorm{\bm{\zeta}}_\infty$ we deduce that
\[
\wabs{N} \geq \wabs{\zeta_k} - \wabs{\phi} - \wnorm{\bm{\psi}}_\infty -
\wabs{\xi} \wlr{\wabs{A} \wabs{\zeta_k} + \wabs{A} \wabs{\phi}
+ \wabs{B} + \wabs{C}}
\]
\[
\geq \wnorm{\bm{\zeta}}_\infty - \wabs{\phi} - \wnorm{\bm{\psi}}_\infty - \frac{0.01}{0.99}
 \wlr{\wnorm{\bm{\zeta}}_\infty + \wabs{\phi} + \wnorm{\bm{\zeta}}_\infty + \wnorm{\bm{\psi}}_\infty}
\]
\[
\geq \wlr{1 - \frac{0.02}{0.99}} \wnorm{\bm{\zeta}}_\infty - \wlr{1 + \frac{0.01}{0.99} }\wlr{\wabs{\phi} + \wnorm{\bm{\psi}}_\infty},
\]
and using the bounds \wref{bound_psi} and \wref{bound_phi} we conclude that
\begin{equation}
\label{lower_n}
\wabs{N} \geq 0.979 \wnorm{\bm{\zeta}}_\infty  - 0.815 \wnorm{\bm{\zeta}}_\infty = 0.164 \wnorm{\bm{\zeta}}_\infty.
\end{equation}
Moreover,  \wref{bound_psi}, \wref{mid_d} and \wref{bound_xia} lead to
\[
\wabs{D - 1} \leq \wabs{\zeta_j} + \wabs{\psi_j} + \wabs{\xi} \wlr{\wabs{A} +\wabs{B} + \wabs{C}}
\]
\[
\leq \wnorm{\bm{\zeta}}_\infty + \wnorm{\bm{\psi}}_\infty + \frac{0.01}{0.99} \wlr{1 + \wnorm{\bm{\zeta}}_\infty + \wnorm{\bm{\psi}}_\infty}
\]
\[
\leq
10^{-3} + 0.403 \times 10^{-3} + \frac{0.01}{0.99} \wlr{1 + 10^{-3} + 0.403 \times 10^{-3} }
\leq 0.012.
\]
Therefore, $ 0.988 \leq D \leq 1.012$.
Combining these bounds on $D$ with \wref{lower_n} we obtain
\[
\wabs{\beta_k} = \wabs{\frac{N}{D}} \geq \frac{0.164 \wnorm{\bm{\zeta}}_\infty}{1.012}
\geq 0.162 \wnorm{\bm{\zeta}}_\infty
\]
and we are done.
\qed{}
\\[0.2cm]

Theorem \ref{thm_lower_bound} is relevant in Salzer's case in Figure \ref{figure_least_squares}
because
\begin{enumerate}
\item[(i)] The columns for $\wnorm{\bm{\zeta}}_\infty$ and $\wnorm{\bm{\zeta}}_\infty/\wlr{\epsilon n}$ in Table 1 below
provide strong empirical evidence that
$2.5 \wlr{n + 3} \epsilon \leq \wnorm{\bm{\zeta}}_\infty \leq 0.001$ for $60 \leq n \leq 1.000.000$ in Salzer's case
in practice.
\item[(ii)] Table 2 in \cite{MascCam}  shows that $\wnorm{\bm{\zeta}}_\infty \leq 0.005$ in this case,
and Lemma 10 in \cite{MascCam} leads to
$S  \leq 1.23 n^2 \leq 5 \times 10^{12}$,
for $S$ in \wref{def_lb_s} and $10 \leq n \leq 2.000.000$. This bound on $S$, the fact that $\wabs{w_j} \geq 1/2$ and equation
\wref{hypo_lb} show that we can apply Theorem \ref{thm_lower_bound}
if $\wabs{x - x_j} \leq 10^{-15}$.
Since $\wabs{x_k} \leq 1$ and $\epsilon \leq 2.3 \times 10^{-16}$, the floating point
number $x \in [-1,1] - \wset{x_j}$ closest to $x_j$ satisfies this condition on $x$.
Therefore, there exist $x_j$ and a floating point number $x$ that satisfies
the hypothesis of Theorem \ref{thm_lower_bound} when $n \leq 10 \leq 2.000.000$.
\item[(iii)] The column for $\wnorm{\bm{\zeta}}_\infty/\wlr{ \epsilon n^2}$ in Table 1 shows that in practice
$\wnorm{\bm{\zeta}}_\infty$ is of order $ \epsilon n^2$ in Salzer's case. This is not surprising because
Lemma 1 in \cite{MascCam} shows that $\zeta_k \approx \sum_{j\neq k} \delta_{jk}$,
the shortest intervals $[x_{k-1},x_k]$ have lengths of order $1/n^2$ and $\hat{x}_k - x_k$ is of order $\epsilon$,
and as a result the largest $\delta_{jk}$ in \wref{def_delta_jk}, and $\wnorm{\bm{\zeta}}_\infty$, are of order $\epsilon n^2$.
\item[(iv)] In summary, Theorem \ref{thm_lower_bound}, in combination with the empirical evidence,
shows that in Salzer's case the maximum backward error for the
barycentric interpolation of Lagrange polynomials grows at least like $\epsilon n^2$, and Theorem \ref{thm_main}
shows that this error grows at most like $\epsilon n^2 \log n$. Therefore, in this case Theorem \ref{thm_main}
is sharp except for a factor of order $\log n$.
\end{enumerate}

We end this section with two tables  presenting the results of experiments with
rounded Chebyshev nodes. The backward errors $\beta$ in these tables are the maximum values
found by evaluating the second barycentric formula in double precision
and comparing the result with the value obtained in
quadruple precision, with $\epsilon \approx 10^{-30}$. For each $n$,
we chose trial points near what we expect to be critical nodes,
as described in Subsection \ref{sub_experiments}.
Table \ref{table_salzer} regards the Salzer's weights $\hat{\wvec{w}}^s$,
and Table \ref{table_rounded} considers the weights obtained by evaluating
numerically $\wfc{\lambda}{\hat{\wvec{x}}^c}$.

\begin{table}[!h]
\caption{The maximum backward error $\beta$ and the relative errors $\bm{\zeta}^s$ in the weights
for Lagrange polynomials with Salzer's weights}
\centering
\begin{tabular}{r|cccccccc}
\hline\\[-0.35cm]
$n$ & $\beta$ & $\wnorm{\bm{\zeta}^s}_\infty$ & $\frac{\beta}{\wnorm{\bm{\zeta}^s}_\infty}$ &
$\frac{\beta}{\epsilon n}$ & $\frac{\wnorm{\bm{\zeta}^s}_\infty}{\epsilon n}$ & $\frac{\beta}{\epsilon n^2}$ &
$\frac{\wnorm{\bm{\zeta}^s}_\infty}{\epsilon n^2}$ \\[0.15cm]
\hline\\[-0.3cm]

$  10 $&        1.9e-15&     8.5e-16&  2.17&     8.3e-01&     3.8e-01&  0.083&  0.038\\
$  20 $&        9.7e-15&     7.1e-15&  1.37&     2.2e+00&     1.6e+00&  0.110&  0.080\\
$  40 $&        1.3e-14&     9.0e-15&  1.45&     1.5e+00&     1.0e+00&  0.037&  0.025\\
$  60 $&        5.9e-14&     5.0e-14&  1.18&     4.4e+00&     3.7e+00&  0.074&  0.062\\
$  80 $&        1.3e-13&     9.0e-14&  1.49&     7.5e+00&     5.1e+00&  0.094&  0.063\\
$  100$&        1.8e-13&     1.6e-13&  1.19&     8.3e+00&     7.0e+00&  0.083&  0.070\\
$  200$&        4.1e-13&     3.1e-13&  1.31&     9.2e+00&     7.1e+00&  0.046&  0.035\\
$  400$&        2.9e-12&     2.3e-12&  1.27&     3.2e+01&     2.6e+01&  0.081&  0.064\\
$  600$&        6.4e-12&     4.1e-12&  1.54&     4.8e+01&     3.1e+01&  0.080&  0.052\\
$  800$&        1.4e-11&     1.2e-11&  1.17&     7.7e+01&     6.5e+01&  0.096&  0.082\\
$  1.000$&      2.4e-11&     2.2e-11&  1.10&     1.1e+02&     9.8e+01&  0.107&  0.098\\
$  2.000$&      3.2e-11&     2.5e-11&  1.28&     7.3e+01&     5.7e+01&  0.036&  0.028\\
$  4.000$&      4.3e-10&     3.8e-10&  1.14&     4.8e+02&     4.2e+02&  0.120&  0.106\\
$  6.000$&      3.5e-10&     1.7e-10&  2.00&     2.6e+02&     1.3e+02&  0.044&  0.022\\
$  8.000$&      1.7e-09&     1.6e-09&  1.09&     9.8e+02&     9.0e+02&  0.122&  0.112\\
$  10.000$&     3.2e-09&     2.7e-09&  1.17&     1.4e+03&     1.2e+03&  0.142&  0.122\\
$  20.000$&     9.4e-09&     7.9e-09&  1.19&     2.1e+03&     1.8e+03&  0.106&  0.089\\
$  40.000$&     6.6e-08&     6.2e-08&  1.06&     7.4e+03&     6.9e+03&  0.185&  0.174\\
$  60.000$&     2.3e-08&     1.7e-08&  1.36&     1.7e+03&     1.2e+03&  0.028&  0.021\\
$  80.000$&     1.3e-07&     1.1e-07&  1.20&     7.5e+03&     6.2e+03&  0.093&  0.078\\
$  100.000$&    1.9e-07&     1.2e-07&  1.63&     8.5e+03&     5.2e+03&  0.085&  0.052\\
$  200.000$&    1.5e-06&     1.1e-06&  1.33&     3.3e+04&     2.5e+04&  0.167&  0.125\\
$  400.000$&    2.6e-06&     1.9e-06&  1.39&     2.9e+04&     2.1e+04&  0.073&  0.053\\
$  600.000$&    4.4e-06&     3.1e-06&  1.40&     3.3e+04&     2.4e+04&  0.055&  0.039\\
$  800.000$&    9.3e-06&     6.9e-06&  1.36&     5.3e+04&     3.9e+04&  0.066&  0.048\\
$  1.000.000$&  8.9e-06&     7.0e-06&  1.26&     4.0e+04&     3.2e+04&  0.040&  0.032
\end{tabular}
\label{table_salzer}
\end{table}

\begin{table}[!h]
\caption{The maximum backward error $\beta$ and the relative errors $\bm{\zeta}^r$ in the weights
for Lagrange polynomials with rounded weights}
\centering
\begin{tabular}{r|cccccccc}
\hline\\[-0.35cm]
$n$ & $\beta$ & $\wnorm{\bm{\zeta}^r}_\infty$ & $\frac{\beta}{\wnorm{\bm{\zeta}^r}_\infty}$ &
$\frac{\beta}{\epsilon n}$ & $\frac{\wnorm{\bm{\zeta}^r}_\infty}{\epsilon n}$ & $\frac{\beta}{\epsilon n^2}$ &
$\frac{\wnorm{\bm{\zeta}^r}_\infty}{\epsilon n^2}$ \\[0.15cm]
\hline\\[-0.3cm]
$  10 $&        9.5e-16&     3.5e-16&   2.70& 0.429& 0.159&     4.3e-02&     1.6e-02\\
$  20 $&        2.4e-15&     6.9e-16&   3.44& 0.535& 0.156&     2.7e-02&     7.8e-03\\
$  40 $&        2.4e-15&     8.3e-16&   2.92& 0.273& 0.094&     6.8e-03&     2.3e-03\\
$  60 $&        3.0e-15&     1.2e-15&   2.55& 0.225& 0.088&     3.8e-03&     1.5e-03\\
$  80 $&        4.8e-15&     2.1e-15&   2.27& 0.270& 0.119&     3.4e-03&     1.5e-03\\
$  100$&        5.2e-15&     2.2e-15&   2.37& 0.235& 0.099&     2.3e-03&     9.9e-04\\
$  200$&        1.2e-14&     5.0e-15&   2.36& 0.263& 0.111&     1.3e-03&     5.6e-04\\
$  400$&        2.0e-14&     9.6e-15&   2.06& 0.223& 0.108&     5.6e-04&     2.7e-04\\
$  600$&        2.5e-14&     1.2e-14&   2.09& 0.190& 0.091&     3.2e-04&     1.5e-04\\
$  800$&        3.6e-14&     1.6e-14&   2.20& 0.202& 0.092&     2.5e-04&     1.1e-04\\
$  1.000$&      4.2e-14&     2.1e-14&   2.02& 0.188& 0.093&     1.9e-04&     9.3e-05\\
$  2.000$&      8.7e-14&     4.3e-14&   2.01& 0.196& 0.097&     9.8e-05&     4.9e-05\\
$  4.000$&      1.7e-13&     8.3e-14&   2.11& 0.196& 0.093&     4.9e-05&     2.3e-05\\
$  6.000$&      2.5e-13&     1.3e-13&   2.02& 0.190& 0.094&     3.2e-05&     1.6e-05\\
$  8.000$&      3.5e-13&     1.6e-13&   2.11& 0.194& 0.092&     2.4e-05&     1.1e-05\\
$  10.000$&     4.1e-13&     2.0e-13&   2.03& 0.187& 0.092&     1.9e-05&     9.2e-06\\
$  20.000$&     8.4e-13&     4.2e-13&   2.00& 0.188& 0.094&     9.4e-06&     4.7e-06\\
$  40.000$&     1.6e-12&     8.2e-13&   2.00& 0.183& 0.092&     4.6e-06&     2.3e-06\\
$  60.000$&     2.5e-12&     1.2e-12&   2.00& 0.184& 0.092&     3.1e-06&     1.5e-06\\
$  80.000$&     3.2e-12&     1.6e-12&   2.00& 0.182& 0.091&     2.3e-06&     1.1e-06\\
$  100.000$&    4.1e-12&     2.1e-12&   2.00& 0.186& 0.093&     1.9e-06&     9.3e-07\\
$  200.000$&    8.1e-12&     4.1e-12&   2.00& 0.183& 0.092&     9.2e-07&     4.6e-07\\
$  400.000$&    1.6e-11&     8.2e-12&   2.00& 0.183& 0.092&     4.6e-07&     2.3e-07\\
$  600.000$&    2.4e-11&     1.2e-11&   2.00& 0.183& 0.091&     3.0e-07&     1.5e-07\\
$  800.000$&    3.3e-11&     1.6e-11&   2.00& 0.183& 0.092&     2.3e-07&     1.1e-07\\
$  1.000.000$&  4.1e-11&     2.0e-11&   2.01& 0.184& 0.092&     1.8e-07&     9.2e-08
\end{tabular}
\label{table_rounded}
\end{table}

Tables \ref{table_salzer} and \ref{table_rounded}, the fact that $\epsilon \approx 2.3 \times 10^{-16}$,
 and the least squares lines in Figure \ref{figure_least_squares}, support the
back-of-the-envelope estimates
\[
\beta_s \approx 0.1 \epsilon n^2
\hspace{1cm} \wrm{and} \hspace{1cm}
\beta_r \approx 0.2 \epsilon n
\]
for the maximum backward errors in the family of Lagrange polynomials,
and Tables \ref{table_salzer} and \ref{table_rounded} show also that $\wnorm{\bm{\zeta}}_\infty$ gives a quite good
estimate of the order of magnitude of the backward error: the ratio $\beta / \wnorm{\bm{\zeta}}_\infty$
is in the range $[1,4]$ for all $n$ considered in both tables
(and we obtained similar results in  literally more than a thousand other experiments.)
Therefore, for large $n$ the $\bm{\zeta}$ for Salzer's weights are considerably larger than
the $\bm{\zeta}$ for rounded nodes.

\subsection{Experimental settings}
\label{sub_experiments}
Our experiments used C++11 code, compiled with g++4.8.1, with usual
options for optimization in release builds: -mavx -O3 -DNDEBUG. We did not
use any tricks to improve performance or accuracy. The experiments were performed
in standard processors: an Intel Core i7-2700K  and an Intel Xeon E5-2640.
The quadruple precision computations were performed with g++'s \_\_float128 type,
which has machine precision of order $10^{-30}$
(we compiled the code with option -fext-numeric-literals  and linked the library quadmath in order to
use these floating point numbers and the constant $\pi$ with precision of $10^{-30}$.)
We checked the results by comparing them with the ones obtained using the MPFR library
(Fousse et al. , 2007)
There were differences in the results obtained by the two
processors. For instance, some nodes computed by the Core i7 have error or order $10^{-16}$
while the same nodes computed by the Xeon have error of order $10^{-18}$, and vice versa.
As a result, some numbers in Tables \ref{table_salzer} and \ref{table_rounded} computed by
these two processors had differences even in their leading digit.
However, all entries for Tables
\ref{table_salzer} and \ref{table_rounded} computed by both processors had the same order of magnitude,
and these tables present the same overall picture in both cases.

For each $n$ in Tables \ref{table_salzer} and \ref{table_rounded} we computed the weights $\hat{\wvec{w}}^r$ and
$\hat{\wvec{w}}^s$ in quadruple precision and then obtained $\bm{\zeta}^r$ and $\bm{\zeta}^s$.
We then choose a set of pairs of indexes $(k,j)$ for Tables \ref{table_salzer} and \ref{table_rounded} as follows:
\begin{itemize}
\item We formed a vector \verb indexes  containing the indexes $0$, $n/2$ and $n$,
the ten indexes corresponding to the largest $z^r_k$, the ten indexes corresponding to the largest
$z^s_k$,  the ten indexes corresponding to the smallest $z^r_k$ and the ten indexes
corresponding to the smallest $z^s_k$, and removed the repetitions.

\item We then formed all the pairs $(k,j)$ with distinct $j$ and $k$ in \verb indexes ,
with $j \neq n / 2$, because $x$ very close to $x_{n/2} = 0$ may lead to underflow.
\item For each pair $(k,j)$ in which $j > 0$, we considered the $5000$ floating point numbers $x$ to the
left of $x_j$. When $j < n$, we considered the $5000$ floating point numbers to the
right of $x_j$. For each trial point $x$ we computed the backward error for the $k$th Lagrange polynomial
evaluated at $x$, for both sets of weights,
and Tables \ref{table_salzer} and \ref{table_rounded} report the maximum backward error found
for each $n$.
\end{itemize}

\section{The stability of the linear map $\winter{x}{w}$}
\label{section_proofs}
The main result in this section is Theorem \ref{thm_delta}, which gives general
bounds on the effects of the perturbations of the nodes and weights of the
second barycentric formula, and leads to the proof of Theorem \ref{thm_main}
presented at the end of the section. In order to allow for general perturbations in the
nodes and in the endpoints of the interpolation interval, this theorem is stated
in terms of a generic map $\chi: [\hat{x}^+,\hat{x}^-] \rightarrow [x^+,x^-]$,
so that the readers could consider maps that would allow them to use
Theorem \ref{thm_delta} in situations which are not considered here. For instance,
it is possible to find a map $\chi$ which allows one to handle cases in which
$x_k = x^-$ and $\hat{x}_k < \hat{x}^-$, but we do not consider such
cases here for the sake of brevity. Lemma \ref{lem_chi} and
its Corollary \ref{cor_chi} present a canonical map $\chi$, which will be appropriate
in most practical situations. The constant $d$ in Theorem \ref{thm_delta}
for this map is the $\delta$ in \wref{def_delta}, and the readers which
are not concerned with utmost
generality can go directly to Theorem \ref{thm_delta}, ignore the
function $\chi$ in the statement of this theorem and take $d = \delta$ in \wref{def_delta}.

\begin{lemma}
\label{lem_chi}
Given numbers $\hat{x}_0 < \hat{x}_1 < \dots < \hat{x}_n$  and $x_0 < x_1 < \dots < x_n$,
the piecewise linear map $\chi: [\hat{x}_0,\hat{x}_n] \rightarrow [x_0,x_n]$ given by
$\wfc{\chi}{\hat{x}_0} := x_0$ and
\[
\wfc{\chi}{\hat{x}} := x_k + \wlr{\hat{x} - \hat{x}_k} \frac{x_{k + 1} - x_k}{\hat{x}_{k+1} - \hat{x}_k}
\hspace{0.5cm} \wrm{for} \hspace{0.5cm}  \hat{x}_k < \hat{x} \leq \hat{x}_{k + 1},
\]
is strictly increasing, $\wfc{\chi}{\hat{x}_k} = x_k$,
$\wabs{\wfc{\chi}{\hat{x}} - \hat{x}} \leq \wnorm{\hat{\wvec{x}} - \wvec{x}}_\infty$ for $\hat{x} \in [\hat{x}_0,\hat{x}_n]$
and, for $0 \leq j \leq n$ and
\begin{equation}
\label{bound_sigma}
\wabs{\frac{\wfc{\chi}{\hat{x}} - x_j}{\hat{x} - \hat{x}_j} - 1}
\leq \max \wset{\wabs{\delta_{jk}}, \wabs{\delta_{j\wlr{k+1}}}},
\end{equation}
for $\hat{x}_k < \hat{x} < \hat{x}_{k+1}$ and $\delta_{jk}$ in \wref{def_delta}.\\[0.1cm]
\end{lemma}


{\bf Proof of Lemma \ref{lem_chi}.} The function $\wabs{\wfc{\chi}{\hat{x}} - \hat{x}}$ is convex in $[\hat{x}_k,\hat{x}_{k+1}]$
and, therefore, its maximum value in this interval is reached at $\hat{x}_k$ or $\hat{x}_{k+1}$, that is,
\[
\max_{\hat{x} \in [\hat{x}_{k},\hat{x}_{k+1}] }  \wabs{\wfc{\chi}{\hat{x}} - \hat{x}} =
\max\wset{ \wabs{\wfc{\chi}{\hat{x}_k} - \hat{x}_k}, \wabs{\wfc{\chi}{\hat{x}_{k+1}} - \hat{x}_{k+1}}}
\]
\[
= \max\wset{ \wabs{x_k - \hat{x}_k}, \wabs{x_{k+1} - \hat{x}_{k+1}}} \leq \wnorm{\hat{\wvec{x}} - \wvec{x}}_\infty.
\]
In order to prove \wref{bound_sigma} it suffices to show that the functions
\[
\wfc{h_{jk}}{\hat{x}} := \frac{\wfc{\chi}{\hat{x}} - x_j}{\hat{x} - \hat{x}_j} - 1 =
\frac{1}{\hat{x} - \hat{x}_j} \wlr{x_k - x_j + \wlr{\hat{x} - \hat{x}_k} \frac{x_{k + 1} - x_k}{\hat{x}_{k+1} - \hat{x}_k}} - 1
\]
satisfy
\begin{equation}
\label{bound_chi_p}
\wabs{\wfc{h_{jk}}{\hat{x}}} \leq \max \wset{ \wabs{\delta_{jk}}, \wabs{\delta_{j\wlr{k+1}}}}
\end{equation}
for $\hat{x}_k < \hat{x} < \hat{x}_{k+1}$. Note that
$\wfc{h_{jk}}{\hat{x}} = A/\wlr{\hat{x} - \hat{x}_j} + B$ for constants
\[
A := x_k  - x_j + \wlr{\hat{x}_j - \hat{x}_k} \frac{x_{k + 1} - x_k}{\hat{x}_{k+1} - \hat{x}_k}
\hspace{0.3cm} \wrm{and} \hspace{0.3cm} B := \frac{x_{k+1} - x_k}{\hat{x}_{k+1} - \hat{x}_k} - 1.
\]
When $j \in \wset{k,k+1}$ we have that $A = 0$, $h_{jk}$ is constant and equal to $\delta_{k\wlr{k+1}}$ in the interval
$\wlr{\hat{x}_{k},\hat{x}_{k+1}}$ and \wref{bound_chi_p} holds. Otherwise, $\hat{x}_j \not \in [\hat{x}_k,\hat{x}_{k+1}]$
and $h_{jk}$ is monotone and continuous in this interval. Therefore, when $j \not \in \wset{k, k+1}$ we
have that
\[
\max_{\hat{x}_k \leq \hat{x} \leq \hat{x}_{k+1}} \wabs{\wfc{h_{jk}}{\hat{x}}} =
\max\wset{ \wabs{\wfc{h_{jk}}{\hat{x}_k}}, \wabs{\wfc{h_{jk}}{\hat{x}_{k+1}}}}
= \max\wset{ \wabs{\delta_{jk}}, \wabs{\delta_{j\wlr{k+1}}}}.
\]
This proves \wref{bound_chi_p} and this proof is complete.
\qed{}
\\[0.2cm]

\begin{corollary}
\label{cor_chi}
Under the conditions \wref{first_cond}--\wref{last_cond},
if $\delta$ in \wref{def_delta}
is  smaller than one then $\hat{x}_k < \hat{x}_{k+1}$ for $0 \leq k \leq n$ and there exists a bijection
$\chi: [\hat{x}^-,\hat{x}^+]   \rightarrow  [x^-,x^+]$
such that $\wfc{\chi}{\hat{x}^-} = x^-$,
 $\wfc{\chi}{\hat{x}_k} = x_k$, $\wfc{\chi}{\hat{x}^+} = x^+$,
 \[
 \wabs{\wfc{\chi}{\hat{x}} - \hat{x}} \leq
 \max \wset{ \wnorm{\wvec{x} - \hat{\wvec{x}}}_\infty, \ \wabs{x^- - \hat{x}^-},  \ \wabs{x^+ - \hat{x}^+}},
 \]
 and, for $0 \leq j \leq n$,
\begin{itemize}
\item If $\hat{x}^- < \hat{x} < \hat{x}_{k^-}$ then
\[
\wabs{\frac{\wfc{\chi}{\hat{x}} - x_j}{\hat{x} - \hat{x}_j} - 1} \leq
\max \wset{ \wabs{\delta^-_j}, \ \wabs{\delta_{jk^-}}}.
\]
\item If $k^- \leq k <  k^+$ and $\hat{x}_k < \hat{x} < \hat{x}_{k + 1}$ then
\[
\wabs{\frac{\wfc{\chi}{\hat{x}} - x_j}{\hat{x} - \hat{x}_j} - 1} \leq
\max \wset{\wabs{\delta_{jk}}, \ \wabs{\delta_{j\wlr{k+1}}}}.
\]
\item If $\hat{x}_{ k^+} < \hat{x} < \hat{x}^+$ then
\[
\wabs{\frac{\wfc{\chi}{\hat{x}} - x_j}{\hat{x} - \hat{x}_j} - 1} \leq
\max \wset{\wabs{\delta_{j k^+}}, \ \wabs{\delta^+_j}}.\\[0.2cm]
\]
\end{itemize}
\end{corollary}

{\bf Proof of Corollary \ref{cor_chi}.}
Corollary \ref{cor_chi} follows from Lemma \ref{lem_chi} applied to the vectors obtained by inserting
$\hat{x}^-$, $\hat{x}^+$,$x^-$, $x^+$ in the appropriate positions of $\hat{\wvec{x}}$ and $\wvec{x}$.
In fact, consider the vectors
\[
\hat{\wvec{x}}' := \wlr{\hat{x}_0,\dots, \hat{x}_{k^- - 1}, \wlrc{\hat{x}^-, \hat{x}_{k^-}}, \hat{x}_{k^- + 1},
\dots , \hat{x}_{k^+ - 1}, \wlrc{\hat{x}_{k^+},\hat{x}^+}, \dots, \hat{x}_{n}}^t,
\]
\[
\wvec{x}' := \wlr{x_0,\dots, x_{k^- - 1}, \wlrc{x^-, x_{k^-}}, x_{k^- + 1},
\dots , x_{k^+ - 1}, \wlrc{x_{k^+},x^+}, \dots, x_{n}}^t,
\]
where $\wlrc{x^-,x_{k^-}}$ represents $x_{k^-}$ when $x^- = x_{k^-}$ and the pair
$x^-, x_{k^-}$ when
$x^- \neq x_{k^-}$, and the other braces are analogous. The hypothesis
$\hat{x}^- = \hat{x}_k \Leftrightarrow x^- = x_k$ and
$\hat{x}^+ = \hat{x}_k \Leftrightarrow x^- = x_k$ ensures that $\hat{\wvec{x}}'$ and $\wvec{x}'$ have the
same dimension, and the definitions of $k^-$ and $k^+$ guarantee that $x'_k < x'_{k+1}$ for the relevant $k$.
Finally, the vector $\hat{\wvec{x}}'$ is strictly sorted because, for instance for $j > k^-$,
\[
\hat{x}_j - \hat{x}^- = \frac{x_{j} - x^-}{1 + \delta^-_j} \geq
\frac{x_{j} - x^-}{1 - \wabs{\delta^-_j}} > 0,
\]
and we can indeed derive Corollary \ref{cor_chi} from Lemma \ref{lem_chi}.
\qed{}\\[0.2cm]

%
%

\begin{theorem}
\label{thm_delta}
Under the conditions \wref{hypo_nz_den} and  \wref{first_cond}--\wref{last_cond},
if $\hat{x} \in [\hat{x}^-,\hat{x}^+] \setminus \wset{\hat{x}_0,\dots,\hat{x}_n}$,
$d \in \wrone{}$ and the function
\[
\chi: [\hat{x}^-,\hat{x}^+] \setminus \wset{\hat{x}_0,\dots,\hat{x}_n} \rightarrow
[x^-,x^+] \setminus \wset{x_0,\dots,x_n},
\]
are such that
\begin{equation}
\label{hypo_main}
\max_{0 \leq k \leq n} \wabs{\frac{\wfc{\chi}{\hat{x}} - x_k}{\hat{x} - \hat{x}_k} - 1} \leq d
< \frac{1 - \wnorm{\wfc{\bm{\zeta}}{\wvec{w},\hat{\wvec{w}}}}_\infty}{\Lambda_{x^-,x^+,\wvec{x},\wvec{w}}} - \wnorm{\wfc{\bm{\zeta}}{\wvec{w},\hat{\wvec{w}}}}_\infty
\end{equation}
 then
\begin{equation}
\label{well_def_s}
\sum_{k = 0}^n \frac{\hat{w}_k}{\hat{x} - \hat{x}_k} \neq 0
\end{equation}
and there exists $\bm{\beta} \in \wrn{n+1}$  such that
\begin{equation}
\label{bound_lem_delta_x}
\wnorm{\bm{\beta}}_\infty \leq
\frac{\wlr{d + \wnorm{\wfc{\bm{\zeta}}{\wvec{w},\hat{\wvec{w}}}}_\infty } \wlr{1 + \Lambda_{x^-,x^+,\wvec{x},\wvec{w}}}}
      {1 - \wnorm{\wfc{\bm{\zeta}}{\wvec{w},\hat{\wvec{w}}}}_\infty - \wlr{d + \wnorm{\wfc{\bm{\zeta}}{\wvec{w},\hat{\wvec{w}}}}_\infty } \Lambda_{x^-,x^+,\wvec{x},\wvec{w}}},
\end{equation}
and
\begin{equation}
\label{back_dx}
\wfc{\winterf{\hat{x}}{\hat{w}}{y}}{\hat{x}} = \wfc{\winterf{x}{w}{\tilde{y}}}{\wfc{\chi}{\hat{x}}}
\hspace{1cm} \wrm{for} \hspace{1cm}
\tilde{y}_k = y_k\wlr{1 + \beta_k}.
\end{equation}
Moreover, if \wref{hypo_main} holds for all $\hat{x} \in [\hat{x}^-,\hat{x}^+]$ then
\begin{equation}
\label{bound_leb_delta_x}
\Lambda_{\hat{x}^-,\hat{x}^+,\hat{\wvec{x}},\hat{\wvec{w}}} \leq
\frac{\wlr{1 + d}\Lambda_{x^-,x^+,\wvec{x},\wvec{w}} }{1 - \wnorm{\wfc{\bm{\zeta}}{\wvec{w},\hat{\wvec{w}}}}_\infty - \wlr{d + \wnorm{\wfc{\bm{\zeta}}{\wvec{w},\hat{\wvec{w}}}}_\infty }
 \Lambda_{x^-,x^+,\wvec{x},\wvec{w}}}.\\[0.3cm]
\end{equation}
\end{theorem}

%
%

{\bf Proof of Theorem \ref{thm_delta}.} Equation \wref{hypo_main} shows that
\[
\nu_k := \frac{\wfc{\chi}{\hat{x}} - x_k}{\hat{x} - \hat{x}_k} - 1
\]
satisfies $\wabs{\nu_k} \leq  d$,
and also that $\wnorm{\bm{\zeta}}_\infty < 1$. Therefore $1 + \zeta_k \neq 0$
and we can write $\hat{w}_k = w_k / \wlr{1 + \zeta_k}$  and deduce that
\[
\sum_{k = 0}^n \frac{\hat{w}_k}{\hat{x} - \hat{x}_k} = \sum\limits_{k = 0}^{n} \frac{w_k}{\wfc{\chi}{\hat{x}} - x_k}
\frac{1}{1 + \zeta_k} \frac{\wfc{\chi}{\hat{x}} - x_k}{\hat{x} - \hat{x}_k}
= \sum\limits_{k = 0}^{n} \frac{w_k}{\wfc{\chi}{\hat{x}} - x_k} \frac{1 + \nu_k}{1 + \zeta_k}
\]
\begin{equation}
\label{den_xy}
= D \wlr{1 + \frac{1}{D} \sum_{k = 0}^n \frac{w_k}{\wfc{\chi}{\hat{x}} - x_k} \frac{\nu_k - \zeta_k}{1 + \zeta_k}}
= D \wlr{1 + E},
\end{equation}
for
\[
D := \sum_{k = 0}^n \frac{w_k}{\wfc{\chi}{\hat{x}} - x_k},
\hspace{0.6cm}
\sigma_k := \frac{\nu_k - \zeta_k}{1 + \zeta_k}
\hspace{0.6cm} \wrm{and} \hspace{0.6cm}
E := \frac{1}{D} \sum_{k = 0}^n \frac{w_k \sigma_k}{\wfc{\chi}{\hat{x}} - x_k}
\]
(Note that \wref{hypo_nz_den} implies that $D \neq 0$.) The bound
$\wabs{\nu_k} \leq  d$ yields
\[
\wabs{\sigma_k} \leq  \frac{d + \wnorm{\bm{\zeta}}_\infty}{1 - \wnorm{\bm{\zeta}}_\infty},
\]
and \wref{hypo_main} and the definition of the Lebesgue constant \wref{def_lebesgue} lead to
\[
\wabs{E} =
\wabs{\frac{1}{D} \sum_{k = 0}^n \frac{w_k \sigma_k}{\wfc{\chi}{\hat{x}} - x_k}} =
\wabs{\wfc{\winterf{x}{w}{\bm{\sigma}}}{\wfc{\chi}{\hat{x}}}}
\]
\begin{equation}
\label{bound_e}
\leq
\Lambda_{x^-,x^+,\wvec{x},\wvec{w}} \wnorm{\bm{\sigma}}_\infty
\leq \frac{d + \wnorm{\bm{\zeta}}_\infty}{1 - \wnorm{\bm{\zeta}}_\infty} \Lambda_{x^-,x^+,\wvec{x},\wvec{w}}  < 1,
\end{equation}
and this bound, in combination with \wref{den_xy}, proves \wref{well_def_s}.
Therefore, $\wfc{\winterf{\hat{x}}{\hat{w}}{y}}{\hat{x}}$ is well defined and
$\wfc{\winterf{\hat{x}}{\hat{w}}{y}}{\hat{x}} = N / (D \wlr{1 + E})$,
for $D$ and $E$ as above and
\[
N :=
\sum_{k = 0}^n \frac{\hat{w}_k y_k}{\hat{x} - \hat{x}_k} = \sum\limits_{k = 0}^{n}
\frac{w_k}{\wfc{\chi}{\hat{x}} - x_k}  y_k \frac{1}{1 + \zeta_k} \frac{\wfc{\chi}{\hat{x}} - x_k}{\hat{x} - \hat{x}_k}
= \sum\limits_{k = 0}^{n} \frac{w_k}{\wfc{\chi}{\hat{x}} - x_k}  \theta_k,
\]
for
\[
\theta_k := \frac{1 + \nu_k}{1 + \zeta_k} y_k .
\]
It follows that
$\wfc{\winterf{\hat{x}}{\hat{w}}{y}}{\hat{x}} = \wfc{\winterf{x}{w}{\tilde{y}}}{\wfc{\chi}{\hat{x}}}$,
with
\[
\tilde{y}_k = \frac{\theta_k}{1 + E} = y_k \wlr{1 + \beta_k}
\hspace{1cm} \wrm{and} \hspace{1cm}
\beta_k := \frac{\frac{1 + \nu_k}{1 + \zeta_k} - 1 - E}{1 + E} = \frac{\sigma_k - E}{1 + E},
\]
and \wref{bound_e} leads to
\[
\wabs{\beta_k} \leq \frac{\wabs{\sigma_k} + \wabs{E}}{1 - \wabs{E}} \leq
\frac{\wlr{d + \wnorm{\bm{\zeta}}_\infty } \wlr{1 + \Lambda_{x^-,x^+,\wvec{x},\wvec{w}}}}
{1 - \wnorm{\bm{\zeta}}_\infty - \wlr{d + \wnorm{\bm{\zeta}}_\infty } \Lambda_{x^-,x^+,\wvec{x},\wvec{w}}},
\]
and we have verified \wref{bound_lem_delta_x} and \wref{back_dx}.
Let us now prove \wref{bound_leb_delta_x}. For each $\hat{x} \in [\hat{x}^-,\hat{x}^+]$, the hypothesis about $d$ and
equation \wref{well_def_s} guarantee that $\wfc{\winterf{\hat{x}}{\hat{w}}{y}}{\hat{x}}$ is well defined
and
\[
\wabs{\wfc{\winterf{\hat{x}}{\hat{w}}{y}}{\hat{x}}} = \wabs{\wfc{\winterf{x}{w}{\tilde{y}}}{\wfc{\chi}{\hat{x}}}}
\leq \Lambda_{x^-,x^+,\wvec{x},\wvec{w}} \wnorm{\tilde{\wvec{y}}}_\infty  \leq
\Lambda_{x^-, x^+,\wvec{x},\wvec{w}} \wnorm{\wvec{y}}_\infty \wlr{1 + \wnorm{\bm{\beta}}_\infty},
\]
and the bound \wref{bound_lem_delta_x} yields
\[
\wabs{\wfc{\winterf{\hat{x}}{\hat{w}}{y}}{\hat{x}}} \leq
\frac{1 + d}
      {1 - \wnorm{\bm{\zeta}}_\infty - \wlr{d + \wnorm{\bm{\zeta}}_\infty} \Lambda_{x^-,x^+,\wvec{x},\wvec{w}}}
\Lambda_{x^-,x^+,\wvec{x},\wvec{w}} \wnorm{\wvec{y}}_\infty.
\]
Taking the sup in $\hat{x} \in [\hat{x}^-,\hat{x}^+]$ of this expression we deduce \wref{bound_leb_delta_x}.
\qed{}\\[0.2cm]

{\bf Proof of Theorem \ref{thm_main}.}
We use Stewart's error counter (Higham, 2002)
\[
\wst{k}  \ \ = \ \  \prod\limits_{i = 1}^{k}(1+\xi_i)^{\sigma_i}, \  \ \ \wrm{for} \ \  \sigma_i \ =  \ \pm 1, \ \ \wrm{and} \ \ |\xi_i| \ \leq \ \epsilon,
\]
and write $\wst{k}_{\ell}$ to give a label $\ell$ to the specific $k$ rounding errors
we are concerned with. We note that \cite{Hig}'s Lemma 3.1 implies that
\[
\wrm{if} \ k \epsilon < 1 \ \ \wrm{then} \ \ \wabs{\wst{k} - 1} \leq \frac{k \epsilon}{1 - k \epsilon}
\hspace{0.5cm} \wrm{and} \hspace{0.5cm} 0 < \wst{k} \leq \frac{1}{1 - k \epsilon}.
\]
If $\hat{x} = \hat{x}_k$ for some $k$ then we can simply take $x = x_k$ and $\tilde{y} = y$.
Let us then consider the case $\hat{x} \in [\hat{x}^-,\hat{x}^+] \setminus \wset{\hat{x}_0,\dots,\hat{x}_n}$.
\cite{Hig1} shows that
\begin{equation}
\label{reduce}
\wrounde{\wfc{\winterf{\hat{x}}{\hat{w}}{y}}{\hat{x}}}
 \ \ =
\ \ \frac{ \sum\limits_{k = 0}^{n} \frac{ \hat{w}_k y_k \wst{n + 3}_{k}}{\hat{x} - \hat{x}_k} }
{\sum\limits_{k = 0}^{n} \frac{\hat{w}_k \wst{n + 2}_{k}}{\hat{x} - \hat{x}_k} } \ \ =
\ \ \frac{ \sum\limits_{k = 0}^{n} \frac{ w'_k y'_k}{\hat{x} - \hat{x}_k} }
{\sum\limits_{k = 0}^{n} \frac{w'_k}{\hat{x} - \hat{x}_k} } \ \ = \ \
\wfc{\winterf{\hat{x}}{w'}{y'}}{\hat{x}},
\end{equation}
for
\[
w'_k := \hat{w}_k \wst{n + 2}_{k} \hspace{1cm} \wrm{and} \hspace{1cm}
y'_k := y_k \frac{\wst{n + 3}_{k}}{\wst{n + 2}_{k}} = y_k\wst{2 n + 5}_{k}.
\]
Recalling that $w_k = \hat{w}_k \wlr{1 + \zeta_k}$ and using Higham's Lemma 3.1 we obtain
\[
\wabs{\zeta_k'} = \wabs{\frac{w_k - w'_k}{w'_k}} =
\wabs{\frac{1 + \zeta_k  - \wst{n+2}_k}{\wst{n + 2}_k}} =  \wabs{\wlr{1 + \zeta_k} \wst{n + 2}_{k'} - 1}
\]
\[
\leq \wabs{\zeta_k \wst{n + 2}_{k'}} +  \wabs{\wst{n + 2}_{k'} - 1} \leq \frac{\wabs{\zeta_k} + \wlr{n + 2} \epsilon}{1 - \wlr{n+2} \epsilon} \leq Z.
\]
Corollary \ref{lem_chi} implies that there exists a function $\chi$ as required
by the hypothesis of Theorem \ref{thm_delta} with $d = \delta$ in \wref{def_delta}, and this
theorem applied to $\hat{\wvec{x}} = \hat{\wvec{x}}$, $\hat{\wvec{w}} = \wvec{w}'$,  and $\wvec{y} = \wvec{y}'$ and equation \wref{reduce} show that
$\wrounde{\wfc{\winterf{\hat{x}}{\hat{w}}{y}}{\hat{x}}} = \wfc{\winterf{x}{w}{\tilde{y}}}{\wfc{\chi}{\hat{x}}}$,
with
\[
\tilde{y}_k := y'_k \wlr{1 + \alpha_k} = y_k \wst{2 n + 5}_k \wlr{1 + \alpha_k},
\]
and
\[
\wabs{\alpha_k} \leq \frac{\wlr{\delta + Z} \wlr{1 + \Lambda_{x^-,x^+,\wvec{x},\wvec{w}}}}
                         {1 - Z - \wlr{\delta + Z} \Lambda_{x^-,x^+,\wvec{x},\wvec{w}}}.
\]
Higham's Lemma shows that $\nu_k := \wst{2 n + 5}_k - 1$ satisfies
\[
\wabs{\nu_k} \leq \wlr{2 n + 5}\epsilon / \wlr{1 - \wlr{2 n+ 5} \epsilon},
\]
and this completes the proof of Theorem \ref{thm_main}.
\qed{}
\bibliographystyle{elsarticle-num}

\end{document}